\documentclass[conference, a4paper]{IEEEtran}
\IEEEoverridecommandlockouts

\usepackage{amsmath,amssymb,amsfonts}
\usepackage{algorithmic}
\usepackage{graphicx}
\usepackage[dvipsnames]{xcolor}
\usepackage[ruled,linesnumbered,noend]{algorithm2e}
\usepackage{breakurl}

\title{A Machine Learning Approach to Boost the Vehicle-2-Grid Scheduling}

\author{
    \IEEEauthorblockN{Gabriele Agliardi\IEEEauthorrefmark{2}\IEEEauthorrefmark{4}, Giorgio Cortiana\IEEEauthorrefmark{1}, Anton Dekusar\IEEEauthorrefmark{2}\IEEEauthorrefmark{3}, Kumar Ghosh\IEEEauthorrefmark{1}, Naeimeh Mohseni\IEEEauthorrefmark{1}, \\ Corey O'Meara\IEEEauthorrefmark{1}, V\'ictor Valls\IEEEauthorrefmark{3}, Kavitha Yogaraj\IEEEauthorrefmark{2}\IEEEauthorrefmark{4}, Sergiy Zhuk\IEEEauthorrefmark{2}\IEEEauthorrefmark{3}}
    \IEEEauthorblockA{
    \IEEEauthorrefmark{1}E.ON Digital Technology GmbH, Essen, Germany,}  
    {\IEEEauthorrefmark{2}IBM Quantum,
    \IEEEauthorrefmark{3}IBM Research Europe -- Dublin, \IEEEauthorrefmark{4}IBM Italy, 
    \IEEEauthorrefmark{5}IBM Research -- India }
}

\begin{document}
\maketitle

%%%%%%%%%% Abstract %%%%%%%%%%
\begin{abstract}
Electric Vehicles (EVs) are emerging as battery energy storage systems (BESSs) of increasing importance for different power grid services. However, the unique characteristics of EVs makes them more difficult to operate than \emph{dedicated} BESSs. In this work, we apply a data-driven learning approach to leverage EVs as a BESS to provide capacity-related services to the grid. The approach uses machine learning to predict how to charge and discharge EVs while satisfying their operational constraints. As a paradigm application, we use flexible energy commercialization in the wholesale markets, but the approach can be applied to a broader range of capacity-related grid services. 
We evaluate the proposed approach numerically and show that when the number of EVs is large, we can obtain comparable objective values to CPLEX and approximate dynamic programming, but with shorter run times. These reduced run times are important because they allow us to (re)optimize frequently to adapt to the time-varying system conditions. 
\end{abstract}

\section{Introduction}

The energy sector has historically relied on centralized and carbon-intense energy sources. However, this is changing with the renewable energy transition. Renewable energy sources are intermittent and probabilistic in nature, and their energy production often peaks when the energy demand is low. This volatility in energy production requires the use of expensive energy sources such as fossil fuels to balance the energy supply, which makes energy prices fluctuate heavily. As a result, the energy industry is investing in battery energy storage systems (BESSs) \cite{iberdrola, flexpwr, en16186546} that can store energy when it is cheap in the wholesale market and sell it back to the grid when energy prices are high.

EVs are emerging as a storage solution of increasing importance. In 2022, the number of EVs on the road exceeded 26 million \cite[pp.\ 14]{ev-outlook-2023}, and the global fleet of EVs is expected to grow to about 240 million in 2030 --- 10\% of the road vehicle fleet \cite[pp.\ 109]{ev-outlook-2023}.
To put this opportunity in perspective, EV batteries can typically store  20 to 100 kWh and can have an output power of 10 kW when connected to a Type 2 charger. Thus, 100 EVs can deliver 1 MWh, which is enough energy to power up to $90$ households per day on average \cite{bonkers}. However, EVs are a storage asset with unique characteristics compared to standard BESS. EVs are not always plugged in, they consume energy, and their batteries must be sufficiently charged depending on the expected mileage ahead. For example, an EV may require to have at least 40 kWh by 7 am. 

BESSs can be used in various grid services, including bulk energy, ancillary, transmission \& distribution infrastructure, end-user/utility customer, and renewable integration \cite{8247485, akhil2013electricity}. The use of EVs/BESSs for different grid services has been extensively studied in the literature \cite{6303105, ALFAVERH2023108949,  liu2014vehicle, PAN202213023}. 
The works can be divided into two categories, depending on whether they exploit the fast reaction time of batteries to help the grid adapt to sudden changes (e.g., frequency regulation \cite{ALFAVERH2023108949, 6303105, liu2014vehicle}), or leverage the large energy storage capacity of EVs/BESSs to provide capacity-related services \cite{PAN202213023, en15197192}.  The two types of problems often require different solutions due to their different operational time scales \cite{cai2017distributed}.

In this work, we study how to leverage EVs as a BESS to provide capacity-related services to the grid. These problems are often modeled mathematically as deterministic Mixed Integer (Non-)Linear Programs (MILP/MINLP) --- e.g., the Unit Commitment (UC) problems \cite{PAN202213023, AHARWAR2023109671, 10106620}. However, in practice, not all the parameters are known in advance, which requires recomputing the EVs charge and discharge schedules as the system changes (e.g., when the number of EVs connected to the grid changes). Such a task can be arduous when the number of EVs is large and EVs have heterogeneous characteristics and complex operational constraints. Furthermore, the time available to (re)compute solutions in practice is often limited to a few minutes, which makes the solver performance critical. To tackle this problem, we apply a learning approach to reduce the time required to operate a fleet of EVs. The approach consists of two steps. First, we use real-world data (e.g., energy prices, energy demand) to simulate a fleet of EVs and identify ``good'' scheduling policies; i.e., how to charge and discharge EVs to minimize a cost function. Second, we use the generated policies to train a machine learning (ML) model that, given a particular scenario, can predict how to charge and discharge EVs while satisfying multiple operational constraints. As a paradigm application, we use flexible energy commercialization in the wholesale markets, but the approach can be applied to a broader range of capacity-related grid services. This paper makes the following contributions:

\begin{itemize}
\item We study the complexity of leveraging EVs as a BESS (Sec~\ref{sec:problem}). Our results show that the problem can be challenging for out-of-the-box commercial solvers (CPLEX) when there are operational constraints that aim to reduce the number of charging cycles of the batteries.

\item We present a learning approach for leveraging EVs as a BESS solution (Sec.~\ref{sec:learning_approach}). The approach uses real-world data to generate scheduling policies and then trains a ML model that predicts how to charge and discharge EVs for a particular scenario. We generate the scheduling policies with approximate dynamic programming \cite{bertsekas1996neuro} and use kernels to learn how to charge and discharge EVs.

\item We evaluate the scaling of the proposed approach numerically. Our results show that when the number of EVs is large, we can obtain a lower cost than an out-of-the-box commercial solver (CPLEX) and a similar cost to approximate dynamic programming, however, with substantially shorter run times. These reduced run times are important because they allow us to (re)optimize frequently to adapt to the time-varying system conditions. 

\end{itemize}

The rest of the paper is organized as follows. Sec.~\ref{sec:problem} introduces the problem and its complexity depending on different operational constraints. Sec.~\ref{sec:learning_approach} presents the proposed learning approach that combines approximate dynamic programming and kernels. In Sec \ref{sec:numerical_experiments}, we present numerical experiments to illustrate our method and compare it with existing solver-based approaches. Finally, in Sec.~\ref{sec:conclusions}, we conclude.

%%%%%%%%%%%%%% Problem statement %%%%%%%%%%%%%%

\begin{figure}[t!]
\centering
\includegraphics[width=0.89\columnwidth]{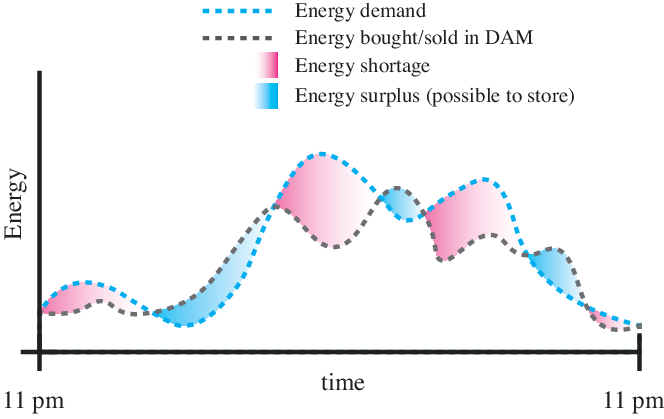}
\caption{Schematic illustration of the problem of leveraging EVs a BESS in the DAM. The blue dashed line indicates the energy demand and the red dashed line is the total energy bought from the markets (DAM + ICM). When the total energy bought is larger than the demand (blue areas), it is possible to store the energy (given that there is enough storage capacity). The red areas indicate an energy shortage that must be covered with storage energy. If there is insufficient energy in storage, the ES must buy more energy from the ICM.}
\label{sec:energy_schema}
\end{figure}

\section{Problem Overview \&  Complexity}
\label{sec:problem}

This section introduces the problem of leveraging EVs as a BESS. We start by giving some background on the energy markets (Sec.~\ref{sec:background}), followed up by two paradigm  business problems from the perspective of energy producers and suppliers (Sec.~\ref{sec:business_problem}). In Sec.~\ref{sec:comp_complexity}, we describe the problem and study its computational complexity for various operational constraints.

%%%
\subsection{Background}
\label{sec:background}

The energy market involves different players, including Energy producers (EPs), Transmission System Operators (TSOs), Distribution System Operators (DSOs), and Energy Suppliers (ESs). EPs sell energy in the wholesale market, which the ESs buy to resell in the retail market---i.e., where consumers (e.g., households and businesses) typically purchase energy. TSOs and DSOs interconnect EPs and ESs, and are responsible for maintaining grid stability.

The wholesale energy market consists of different markets \cite{epexspot}. The Day-Ahead Market (DAM) is a daily auction where exchange members may trade electricity in 24 one-hour trading periods. The Intra Day Auctions (IDA) and Intraday Continuous Market (ICM) are markets that trade energy in 15/30 minute intervals. Typically, the energy in the DAM is cheaper than the energy sold in the IDA and ICM.

\subsection{Trading energy in wholesale markets} 
\label{sec:business_problem}

The problem of trading energy in wholesale markets can be viewed from the perspective of different players. EPs aim to sell as much energy as possible in the DAM market, provided it is possible to deliver it the next day. The challenge they face is that energy production at a given time is uncertain due to the intermittent and probabilistic nature of renewable energy sources. If energy production is lower than the energy sold, EPs need to steer energy production with fossil fuels, energy in storage, or buy energy from other markets (e.g., ICM). 

Similarly, ESs buy electricity in the DAM to meet their customer's demands (i.e., the consumers). The amount of energy bought is based on a forecast of the consumers' demand. ESs always have to deliver the amount of energy the consumers request. Thus, if the demand exceeds the energy bought in the DAM, ESs must obtain energy elsewhere; for example, they can use energy in storage or purchase energy in the wholesale markets (IDAs, ICM). The goal of ESs is to maximize revenue by buying energy from the markets at the lowest possible price (while meeting the consumers' demand). 

In this paper, we consider the following general \emph{operational} problem. Given the energy already bought/sold by an ES/EP in the wholesale market, we want to leverage EVs as a BESS to \emph{minimize the cost of delivering energy.} Fig.~\ref{sec:energy_schema} shows a schematic illustration of the problem. The gray dashed line represents the energy bought/sold in the wholesale market, and the blue line represents the energy demand/energy committed. The blue areas are energy surpluses that can be stored in BESSs, such as in EVs. The red areas are energy shortages that need to be covered either with the energy in storage or by buying energy in the intra-day markets.

\subsection{Computational complexity of leveraging EVs as a BESS}
\label{sec:comp_complexity}

\subsubsection{Setting with increased levels of complexity}
Consider the case where an EP has sold energy in the DAM---in MWh in one-hour intervals. The EP has access to a fleet of EVs to use as a BESS from 6 pm to 6 am the next day. Every 15 minutes, the EP decides (i) how much energy to buy in the ICM (in MWh) and (ii) which EVs should store any excess energy after the committed volume is delivered. The energy stored in the EVs is in kWh, and an EV can charge and discharge between -10 kW and 10 kW. A basic requirement is that the state of charge of an EV must always be within a permissible range; for example, between 0 and the maximum state of charge (kWh) or between 10\% and 80\% of its maximum state of charge. As described in Sec.~\ref{sec:business_problem}, the goal is to minimize the cost of delivering the committed volume (i.e.,  the energy sold in the DAM).

The computational complexity of solving the problem above depends on the characteristics/features we consider for the EV BESS. In the following, we use the methodology described in \cite[Sec.~VII]{abbas2023quantum} and consider scenarios with \emph{increasing} levels of complexity.

\begin{itemize}
\item \textbf{Level 0.} All EVs are plugged in and out at the same time (e.g., a fleet of delivery vehicles), and they all have the same initial state of charge, capacity, and state of charge requirements.

\item \textbf{Level 1.} EVs may have different initial states of charge, desired states of charge, and battery capacities.

\item \textbf{Level 2.} In addition to Level 1, the charging site constrains how much a group of EVs can charge or discharge.

\item \textbf{Level 3.} In addition to Level 2, the EV charging cycles are capped to extend their lifespan.    
\end{itemize}

At level 0, a fleet of EVs is equivalent to having one EV with a ``large battery.'' Level 1 makes all EVs different, increasing the number of variables and constraints required to model and solve the problem. Level 2 couples EVs charging and discharging decisions, breaking the separability of the problem that may be exploited by some solvers. Level 3 adds precedence constraints: the decision of how an EV charges/discharges at a given time affects how they can charge/discharge in the future. This constraint is similar to the min-up-time and min-down-time constraints in UC problems. 

It is possible to add more features or levels to increase the problem's complexity. For example, we could consider that EVs can plug in and out instead of always being connected to the grid, or that the electricity demand is not known. However, as we will show next, Level 3 makes the problem already complex to motivate the approach in Sec.~\ref{sec:learning_approach}. 

\begin{figure}
\centering
\input{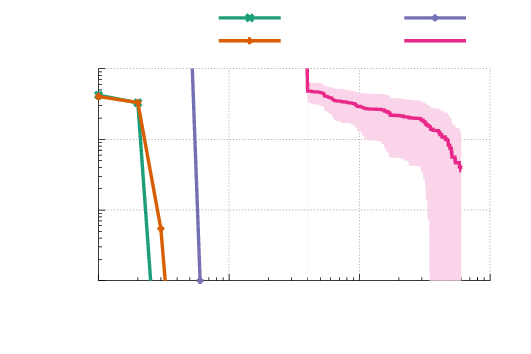}
\caption{CPLEX optimality gap in euros as a function of the run time for 200 EVs and different levels of complexity. The optimality MIP gap is difference between the \emph{best integer solution} and \emph{best (continuous) lower bound}. The results are the average of 50 days from EPEXSPOT \cite{epexspot} selected uniformly at random (2015-2019). The shaded area indicates the standard deviation. The solver is run for 10 minutes and the optimality gap is collected every second.
}
\label{fig:cplex_levels}
\end{figure}

\subsubsection{Problem computation complexity with CPLEX}
\label{sec:complexity}

Consider a scenario with a fleet of 200 EVs that are always connected to the grid from 6 pm to 6 am the next day (i.e., 12 hours). EVs have 100 kWh batteries and can charge/discharge at a rate of $\pm$10 kW. Also, EVs require that their batteries have at least 10 kWh at all times. The optimization model needs to provide the charge/discharge schedule of each vehicle, with a 15-minute resolution. In \textbf{Level 0,} we consider that EVs have an initial state of charge of 30 kWh and must be charged to 80 kWh by 6 am. For \textbf{Level 1,} we select the initial state of charge uniformly at random between 10 and 50 kWh, and the final state between 70 and 100 kWh. For \textbf{Level 2,} we consider that the EVs cannot collectively charge nor discharge more than 1800 kW, namely 90\% of the charging capacity of the fleet. For \textbf{Level 3,} every EV is allowed to switch from charge to discharge (or from discharge to charge) at most once an hour.

We model the problem as a MILP in CPLEX, and show, in Fig.~\ref{fig:cplex_levels},  the evolution of the optimality MIP gap\footnote{The difference between the best integer solution found and the best lower bound coming from the relaxed continuous problem.} (in euros) depending on the runtime. The results are the average of 50 different realizations/days. Observe from the figure that, with our MILP formulation\footnote{The efficiency of CPLEX (or other solvers) depends on the problem formulation.} and with standard solver settings, CPLEX can find an exact solution in a few seconds with Levels 0-2. However, with Level 3, CPLEX needs about 40 seconds to achieve a gap of 4,800 euros (24 euros per EV, corresponding to an 89\% relative gap)\footnote{The relative MIP gap at a given time is the ratio between the MIP gap and the best integer solution found at that time.}, and 10 minutes to drop the gap down to 360 euros (1.8 euros per EV, 13\% relative gap). We believe this is because the constraints introduced at Level 2 and Level 3 make the feasible solutions ``far away'' from each other, which complicates the search in CPLEX.

\section{Learning Approach}
\label{sec:learning_approach}

In this section, we present a learning approach for solving the problem of using EVs as  a BESS. We start by presenting a mathematical model of the problem.

\subsection{Mathematical model}
\label{sec:control_model}

Time is divided in slots of equal duration. We use $k \in \{0,\dots,N-1\}$ to denote the index of a time slot and $N$ the number of times slots we operate the system. The duration of a slot in our problem is typically 15 minutes and the horizon 12 hours (i.e., $N=48$). That is, the duration of a time slot is aligned with the frequency in which we can buy/sell energy in the wholesale market.

We use  $d$ to denote the number of EVs, and $i_k$ to indicate the state of the system at time $k \in \{0,\dots,N-1\}$. The state of the system contains information such as the EVs connected, their state of charge, energy prices, etc. The state of the system  affects the actions that the system controller can make. For example, $i_k$ can be a collection of $d$ pairs indicating the EVs connectivity and their state of charge. 

In each time slot $k \in \{0,\dots,N-1\}$, the controller selects an action $\mu_k$ from the action set 
\[
U(i_k) \subseteq U
\]
where $U$ is a finite set of all possible actions. For example, an element in $U(i_k)$ can be a vector indicating the kWh that an EV can charge or discharge during a time slot. If an EV is not connected in state $i_k$, then $U(i_k)$ will not contain an action to charge or discharge the EV. Similarly, if an EV is fully charged in state $i_k$, $U(i_k)$ will contain the action to discharge the EV but not to charge it.

We want to design a policy $\pi$ that generates a sequence of actions $\{ \mu_0, \mu_1, \dots \}$, $\mu_k \in U(i_k)$ for all $k \in \{0,1\dots, N-1\}$. The  action taken at time $k$ affects the cost and the state in the next time slot. In particular, function $f$ maps a state $i_k$ and an action $\mu_k$ to the next state
\[
i_{k+1} = f(i_k, \mu_k).
\]
Next, let $\rho_k$ denote the price of buying energy from the ICM at time $k \in \{0, \dots, N-1\}$ in euros per kWh. Similarly, let $v_k \in \mathbf R$ be the energy that needs to be delivered at time $k$ in kWh that \emph{cannot} be met with the energy the EP generates at time $k$ (equivalently, an the energy an ES has bought in the DAM). If $v_k < 0$, this indicates that the energy produced is larger than the energy sold. The cost function is given by
\begin{align*}
g(i_k, \mu_k) & := \rho_k \max \{ 0, v_k + \langle \mathbf 1, \mu_k \rangle \}  \\
& \qquad  -   \sigma_k \min \{ 0, v_k + \langle \mathbf 1, \mu_k \rangle \},
\end{align*}
where $\mathbf 1$ is the all-ones vector. Note that $\max \{ 0, v_k + \langle \mathbf 1, \mu_k \rangle \}$ is the energy that needs to be bought in the ICM,\footnote{e.g., when it is not possible to deliver the committed volume} and $ \min \{ 0, v_k + \langle \mathbf 1, \mu_k \rangle \}$ is the excess energy available.\footnote{$\langle x , y \rangle$ denotes the inner product of vectors $x$ and $y$.} The scalar $\sigma_k > 0 $ indicates the cost, in euros, of having excess energy. 
The optimization problem we want to solve is the following:
\begin{align}
\begin{tabular}{ll}
$\underset{u_k \in U(i_k)}{\text{minimize}}$ & $\displaystyle \sum_{k=0}^{N-1} g(i_k, u_k)$ \\
subject to & $i_{N} \in I$,  $i_{k+1} =f(i_k, u_k)$
\end{tabular}
\label{eq:optimization}
\end{align}
where $I$ the set of admissible terminal states. For example, the states that indicate that the EVs are charged enough according to the EV owners requirements. The goal is to minimize the total cost of buying energy from the market, subject to the constraint that the EVs final state is admissible (e.g., all EVs have at least 90 kWh). 

To handle continuously changing inputs, such as energy prices and vehicle availability, the problem in Eq.~\eqref{eq:optimization} is re-optimized every 15 minutes. Unfortunately, 15 minutes may not be enough time to re-optimize and obtain ``good enough'' solutions with commercial solvers, as we showed in Sec.~\ref{sec:complexity}.

%%%%%%%%%%%%

\subsection{Proposed approach: Approximate DP +  Learning}

\begin{figure}
    \centering
    \includegraphics[width=0.8\columnwidth]
    {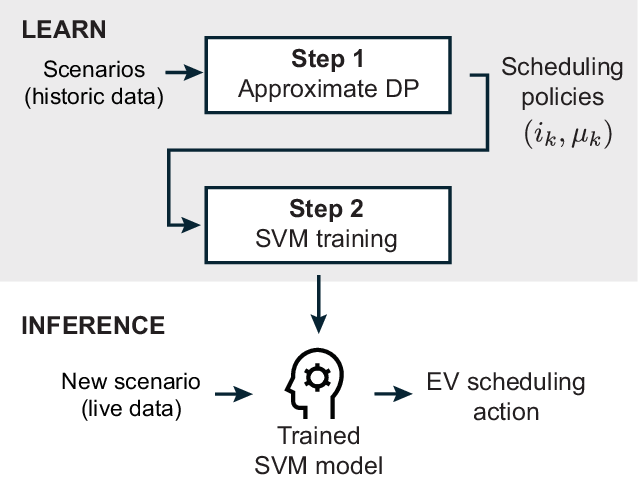}
    \caption{
    Schematic illustration of the proposed approach.}
    \label{fig:workflow}
\end{figure}

To overcome the time limitation mentioned in the previous section, we apply a learning approach that consists of two steps, see Fig.~\ref{fig:workflow}. The first step is to generate scheduling policies by solving the optimization in Eq.~\eqref{eq:optimization} with historical data (energy prices, energy demand, EVs connectivity/charge, etc.). The data generated consists of pairs $(i_k, \mu_k)$ that map the state at time $k$ to an action of charging or discharging EVs. The second step is to train a ML model that, given a state $i_k$, can predict how to charge and discharge EVs to minimize the total cost. %Once a new real-world scenario needs to be optimized, the ML model is applied iteratively over time, starting from the current state $i_0$ to produce the action $\mu_0$ and therefore the subsequent state $i_1$, etc. More precisely, the ML model predicts a decision at the aggregate level, which is then translated into a decision at the vehicle level through a heuristic, as we further detail afterwards in this subsection. 
Notably, there is no severe restriction on the time needed to generate the scheduling policies and to train the ML model, since this is done offline.

The generation of the training data for the ML model (i.e., the scheduling policies) requires to solve a set of input scenarios, and it can, in principle, be performed by multiple optimizers. Here, we propose to use approximate dynamic programming (DP) \cite{bertsekas1996neuro} since it is a well-established framework, strongly connected to branch-and-bound approaches used by commercial solvers. In addition, it also allows us to model stochastic problems—unlike branch-and-bound.
For predicting how EVs should charge/discharge, we propose to use a support vector machines (SVMs). SVMs are appealing because they can make use of kernel functions to efficiently treat non-linearities in data. However, other machine learning approaches are possible, e.g., neural networks.

\subsubsection{Approximate dynamic programming}
\label{sec:dynamic_programming}

We start by describing how we would solve the problem \emph{exactly} with DP, and then show how to solve it approximately.  

Without loss of generality, let us assume that EVs are always connected. Thus, the choice of action at time $k$ determines the state at time $k+1$ deterministically. 
Define the cost-to-go function \cite{bertsekas1996neuro}: 
\begin{align}
J_N^\pi(i) = G(i_N) + \sum_{k=0}^{N-1} g(i_k, \mu_k, i_{k+1}) \quad \text{with $i_0 = i$} 
\label{eq:cost-to-go}
\end{align}
where $\mu_k = \pi(i_k)$ and $G(I_N)$ is the cost of terminating in state $i_N$. In particular, $G(i_N) = +\infty$ if $i_N \in I$ and $0$ otherwise. The optimal $N$-stage cost-to-go starting from state $i$ is 
\begin{align}
J_N^* := \min_{\pi} J _N^\pi (i).
\end{align}
For a given initial state $i_0$, the goal is to minimize $J_N(i_0)$. That is, find
\begin{align}
J_N^*(i_0) = \min_{u \in U(i_0)} \left\{ g(i_0, u, i_1) + J_{N-1}^*(i_1) \right\}
\end{align}
where $J_{N-1}^*(i_1)$ is defined in Eq.\ \eqref{eq:cost-to-go}, and $\mu_0$ in our policy is the solution to the optimization above.

Algorithm \ref{al:DP} contains a DP implementation for the problem. The algorithm finds an \emph{optimal} solution to the optimization problem provided it is feasible. Informally, the program works ``backwards'' to find ``a path'' from the final state $i_N$ to the initial state $i_0$ with minimum cost. Such an approach is different from brute-force search because when this detects a suboptimal strategy, the program stops exploring such a path. Unfortunately, the runtime of Algorithm \ref{al:DP} is typically exponential because of the (i) the number of actions and/or (ii) the number of charging states EVs can be. For instance, if an EV has three different options (discharge, do nothing, and charge), then the number of actions that DP needs to explore for each states is bounded by $2^d \le |U(i)| \le 3^d$ when all EVs are connected. 

To reduce the runtime of DP, we propose to select the action from a subset $\tilde U(i_k)$ of $U(i_k)$.
The solution to the DP may not be the optimal with $\tilde U(i_k)$. However, by choosing  $\tilde U(i_k)$ appropriately, it is possible to obtain  good enough approximations. Specifically, at each time, only three actions are evaluated, namely \texttt{C} (charge as many vehicles as possible, while guaranteeing the existence of a feasible path covering the entire horizon), \texttt{I} (keep idle on all vehicles), or \texttt{D} (discharge as many vehicles as possible, while guaranteeing feasibility). If, say, the policy \texttt{C} is achieved by different choices of \textit{which} vehicles to charge, one choice only is selected based on a heuristic rule. More precisely, in such case our rule would charge the vehicles with the lowest battery. The restriction to the three actions provides satisfactory results, as we demonstrate empirically in the next section.

\begin{algorithm}[t!]
    
    \SetInd{.5em}{.5em}
    \label{al:DP}
    \caption{Dynamic programming pseudo-code}
    \KwData{Cost function; initial and final states of charge of EVs ($i_0$ and $i_N$ resp.); State dynamics function; Constraints.}
    \KwResult{Charge and discharge schedule, i.e., $\{ \mu_0, \mu_1, \dots, \mu_{N-1}\}$}
    %Compute initial state $i_0$ \\
    %Compute final state $i_N$\\
    $J_0(i_N) \leftarrow 0$ ;
    $M_N(i_N) \leftarrow \emptyset$\\
    Store $(N, i_N, J_0(i_N), M_N(i_N))$
    \\
    \For{$k \in \{N,\dots,1\}$}{
    	\For {$j \in (k, j, \cdot, \cdot)$}{
        % Algorithm steps here
        \For{all actions $u \in U$}{
             Compute previous state $i$, i.e., the state $i$ for which it is possible to go to state $j$ with action $u$
            %$j \leftarrow i - j$ \# linear dynamics
            
            \If{state $i$ is feasible according to constraints}{
            $\text{cost} \leftarrow g(i,u,j) + J_{N-k}(j)$
            }
             \If {state $i$ has \emph{not} been seen before at time $k-1$; or state $i$ has been seen before and \textup{cost} is smaller than $J_{N-k+1}(j)$}{
             $J_{N-k+1}(j) \leftarrow \text{cost}$ \\
             $M_{k-1}(i) \leftarrow u \cup M_k(j)$ \\
             Store $(k-1, i, J_{N-k+1}, M_{k-1})$
              }
        }
        }
}
\If{ $(0, i_0, \cdot, \cdot)$ exists}{
\Return $M_0(i_0) = \{\mu_0, \mu_1, \dots, \mu_{N-1} \}$ 
}
\Return Problem is infeasible

\end{algorithm}

\subsubsection{Learning the EV scheduling actions with kernels}

Since a ML model would hardly be able to predict actions at EV granularity, we employ the same restriction as above to the three actions (\texttt C, \texttt I, and \texttt D), and we use the same heuristic as before to determine a corresponding action for each EV. The machine learning problem is thus reduced to a 3-class classification, and we employ kernel methods. Kernel methods are often used to map input data into a higher-dimensional space where this becomes linearly separable. 

\section{Numerical Experiments}
\label{sec:numerical_experiments}

In this section, we evaluate the proposed approach numerically. We consider multiple scenarios with fleets up to 300 EVs and the Level 3 complexity defined in Sec.~\ref{sec:comp_complexity}.

\subsection{Setup}

We generate the training data with approximate DP (Sec.~\ref{sec:dynamic_programming}) with different fleet sizes and price traces from EPEXSPOT \cite{epexspot} (2015-2019).

The prediction model is a Support Vector Classifier (SVC) or its Quantum counterpart (QSVC) \cite{havlicek_supervised_2019}. The SVC or QSVC model is trained with 1000 samples selected uniformly at random. The set of features $i_k$ contains the current time $k$, the state of charge of EVs, and the energy volumes and prices. Volumes are normalized over the number of EVs. Out of those features, only the top 5 most relevant are selected by assigning feature importance with a tree classifier. In terms of labels, three classes are defined, depending the optimal aggregate action obtained by the approximate DP (\texttt{C}, \texttt{I}, or \texttt{D}). The same trained prediction model is used for all number of EVs and for all time periods. The model is then tested on 100 scenarios for each considered fleet size. The test scenarios are chosen from the same pool above, randomly and independently of the training sets. The SVC implementation in scikit-learn and the QSVC in Qiskit are used, with standard parameters. For the QSVC, the default ansatz is adopted, namely a ZZ feature map with full entanglement, 2~repetitions, and one qubit per feature. QSVC is run on the statevector simulator.

\subsection{Results}
\begin{figure}
\centering
\resizebox{0.77\columnwidth}{!}{\input{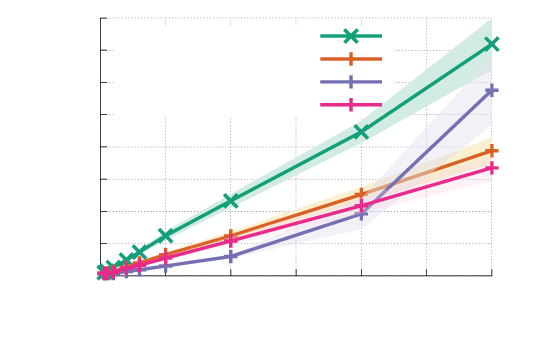}}
\vspace{-0.1cm}
\caption{Value of the objective function of the best solution obtained with different methods (average over 100 scenarios), when the number of EVs varies. Shaded areas indicate 0.2 standard deviation. Lines are drawn for the proposed SVC/QSVC-based algorithm, for the approximate DP used in data generation, and for CPLEX with a runtime limited to 10 minutes. 
}
\label{fig:obj_fn}
\end{figure}

\begin{figure}
\centering
\resizebox{0.77\columnwidth}{!}{\input{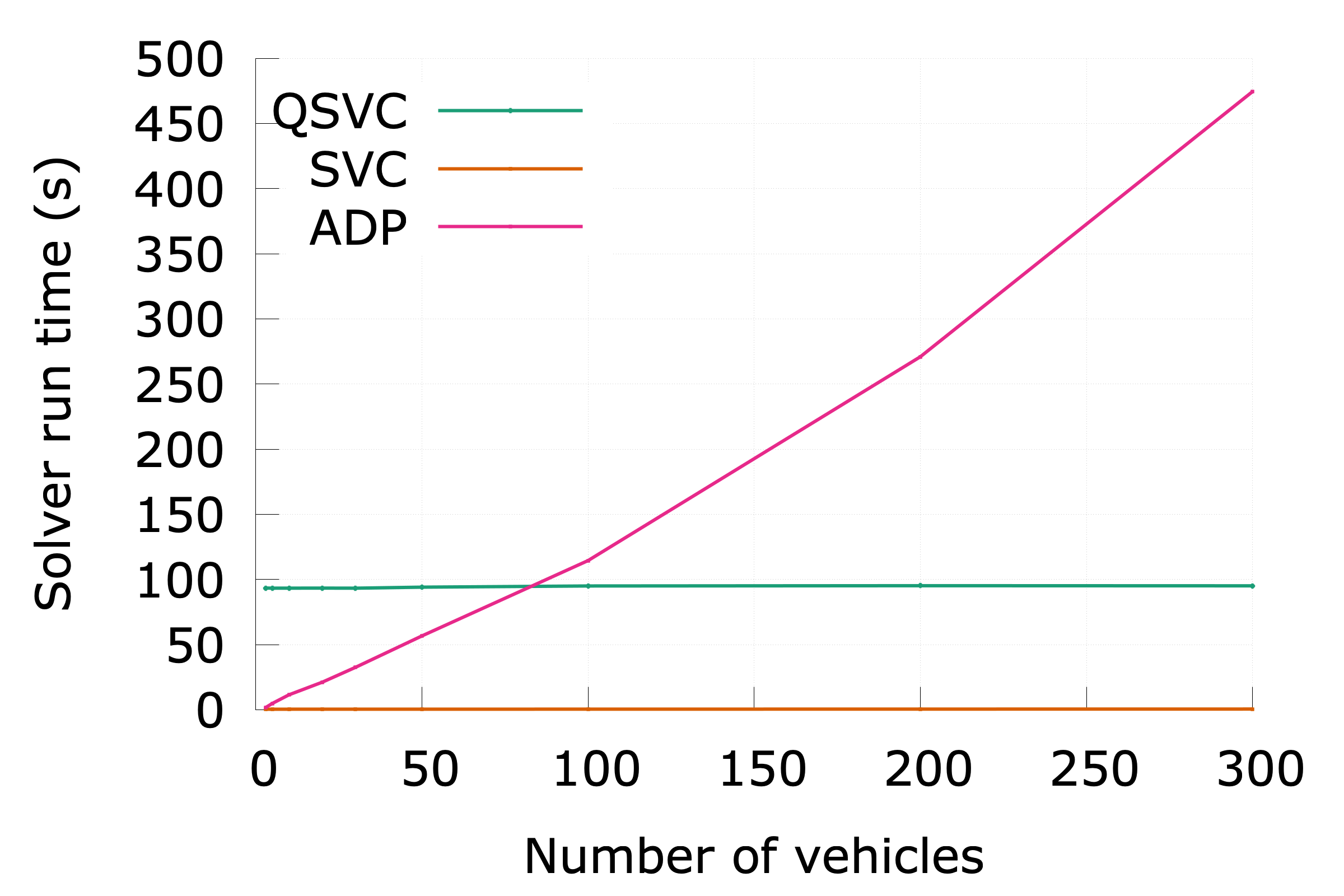}}
\vspace{-0.1cm}
\caption{Run times of the methods in Fig.~\ref{fig:obj_fn}, excluding CPLEX which is time-bound to 10 minutes. Run time is measured as wall-clock time, and is restricted to inference since training is performed offline. Shaded areas indicate the standard deviation, which is inappreciable for SVC and QSVC.
}
\label{fig:time_scaling}
\end{figure}

Fig.~\ref{fig:obj_fn} compares the cost obtained by different methods. Our MILP formulation achieves the best  (i.e., lowest) objective value at up to 200 EVs, when run on CPLEX with standard settings. With 300 EVs, CPLEX hits the 10-minute runtime limit in many problem instances and, as a result, performance degrades heavily. The SVC performs slightly worse than CPLEX for less than 200 EVs, and it outperforms CPLEX when the number of EVs is approximately larger than 230. QSVC obtains worse objective values than SVC due to its classification score being 87\%, compared to 92\% of SVC. Finally, we also show the approximate DP we use to generate data. Approximate DP obtains objective values slightly better than SVC, however, at the price of slower run times. Observe from Fig.~\ref{fig:time_scaling} that the run time of approximate DP scales linearly with the number of EVs, whereas the run time of SVC and QSVC is almost insensitive to the number of EVs. Hence, we use approximate DP only for the training data generation in the learning stage (i.e., the schedules) and not for selecting schedules in ``real-time.'' The run time of QSVC is shown in the figure, despite its simulator origin, with the sole purpose of demonstrating its scalability when the number of EVs increases. 

In summary, the results show that SVC-based optimization obtains good approximate solutions fast. The QSVC-based technique performs worse than its classical counterpart, but it can be faster than approximate DP when the number of EVs is large. 

%%%%%%%%%%% Conclusions %%%%%%%%%%%
\section{Conclusions}
\label{sec:conclusions}

EVs are emerging as a promising BESS solution. However, operating EVs is challenging due to their sophisticated operational constraints and the dynamic nature of the power grid. To address this challenge, we have applied an approach that employs machine learning to predict optimal charging and discharging strategies for EVs while satisfying their multiple operational constraints. Our numerical experiments show that our approach obtains a similar cost as CPLEX and approximate DP but with better run time and scalability. 

%%%%%%%%%%% References %%%%%%%%%%%

\bibliographystyle{IEEEtran}
\bibliography{references}

\end{document}